\documentclass[11pt,english]{article}
\usepackage[T1]{fontenc}
\usepackage[latin1]{inputenc}

\makeatletter

\newcommand{\binom}[2]{{#1 \choose #2}}

\title{The Semigroup of a Word}
\author{Peter M. Higgins}
\date{}
\def\qed{\quad\vrule height4.17pt width4.17pt depth0pt}

\makeatother

\usepackage{babel}
\begin{document}

\title{\textbf{Permutations of a semigroup that map to inverses}}

\author{Peter M. Higgins, University of Essex, U.K.}
\maketitle
\begin{abstract}
We investigate the question as to when the members of a finite regular
semigroup may be permuted in such a way that each member is mapped
to one of its inverses. In general this is not possible. However we
reformulate the problem in terms of a related graph and, using an
application of Hall's Marriage Lemma, we show in particular that the
finite full transformation semigroup does enjoy this property.
\end{abstract}

\section{Introduction and Background }

One of the major subjects of John Howie's research was the full transformation
semigroup {[}11 - 15{]}. The topic is flavoured very much by whether
the base set of the semigroup of mappings is infinite or finite. In
this article we continue the latter thread by studying a question
that may be asked of any semigroup but which is more natural in the
finite regular case and proves particularly interesting for the finite
full transformation semigroup. 

The author's paper {[}6{]} in \emph{Semigroup Forum} was dedicated
to the late John Howie as it represented a continuation of our work
in {[}8,9,10,11{]} and especially {[}10{]} written jointly with John's
colleagues Nik Ruskuc and James Mitchell from St Andrews University
where John worked for most of his career. A feature of the paper was
the application of Hall's Marriage Lemma, a theme that also continues
here.

Throughout, $S$ denotes a semigroup and $T_{X},PT_{X}$ the full
and partial transformation semigroups respectively on the base set
$X$. When $X=X_{n}=\{1,2,\cdots,n\}$ we write $T_{n}$ for $T_{X}$.
Following the texts of Howie {[}16,17{]} and the author's {[}7{]}
we denote the set of idempotents of $S$ by $E(S)$ and the set of
regular elements of $S$ by Reg$(S)$; we shall write $(u,v)\in V(S)$
if $u$ and $v$ are mutual inverses in $S$. We shall also denote
this as $v\in V(u)$ so that $V(u)$ is the set of inverses of $u\in S$.
We extend the notation for inverses to sets $A$: $V(A)=\bigcup_{a\in A}V(a)$.
Standard results on Green's relations, particularly those stemming
from Green's Lemma, will be assumed (Chapter 2 of {[}17{]}, specifically
Lemma 2.2.1) and indeed basic facts and definitions concerning semigroups
that are taken for granted in what follows are all to be found in
Howie's excellent book {[}17{]}. 

Let $C=\{A_{i}\}_{i\in I}$ be any finite family of finite sets (perhaps
with repetition of sets). A set $\tau\subseteq\bigcup A_{i}$ is a
\emph{transversal }of $C$ (or a \emph{system of distinct representatives})
if there exists a bijection $\phi:\tau\rightarrow C$ such that $t\in\phi(t)$
for all $t\in\tau$. We shall denote this by $\tau\,t\,C$. We also
assume Hall's Marriage Lemma in its various forms. In particular $C$
has a transversal if and only if \emph{Hall's Condition }is satisfied,
which says that for all $1\leq k\leq|I|$, the union of any $k$ sets
from $C$ has at least $k$ members (for this and related background
see, for example, the text {[}22{]}).

\textbf{Definitions 1.1 }Let $S$ be any semigroup and let $F=\{f\in PT_{S}:f(a)\in V(a)\,\forall a\in\mbox{dom\,\ \ensuremath{f}\}}$.
We call $F$ the set of \emph{partial inverse matchings }of $S$.
We call $f\in F$ a \emph{permutation matching }if $f$ is a permutation
of $S$; more particularly $f$ is an \emph{involution matching} if
$f^{2}=\iota$, the identity mapping.

We note that, by definition, $S$ is regular if and only if $S$ has
a (full) matching. The word `matching' here does not necessarily have
the connotation of `mutual pairing' as it does in graph theory but
its use in this context is nonetheless convenient. Some other simple
observations are:
\begin{itemize}
\item An inverse semigroup $S$ has a unique matching $a\mapsto a^{-1}$,
which is an involution.
\item A completely regular semigroup $S$ (a semigroup that is a union of
its subgroups) is a disjoint union of its maximal subgroups and so
$a\mapsto a^{-1}$, where $a^{-1}$ is the (unique) group inverse
of $a$, defines an involution matching of $S$. 
\item The identity mapping is a permutation matching if and only if $S$
satisfies the equation $x=x^{3}$; in particular this is true of bands.
\item $S$ has the property that all of its permutations are permutation
matchings if and only if $S$ is a rectangular band as that is exactly
the class of semigroups where every member is inverse to every other
(characterized by either of the equivalent identities $x=xyx$ or
$x=xyz$).
\end{itemize}
Semigroups with an involution matching $(*)$, which satisfies the
additional property that $(ab)^{*}=b^{*}a^{*}$ are called \emph{regular
{*}-semigroups. }They comprise an umbrella class for the classes of
inverse semigroups and completely regular semigroups. There are however
other collections of regular {*}-semigroups of interest including
that of \emph{partition monoids} (see for example {[}3{]}).

In the remainder of the paper we shall, unless otherwise stated, assume
that $S$ is regular and finite. We shall often denote a matching
simply by $'$, so that the image of $a$ is $a'$. The matching $(\cdot)'$
is then permutative if $S'=\{a':a\in S\}=S$ and $(')$ is an involution
matching if $(a')'=a\,\forall a\in S$. 

We shall work with the family of subsets of $S$ given by $V=\{V(a)\}_{a\in S}$.
The members of $V$ may have repeated elements\textemdash for example
$S$ is a rectangular band if and only if $V(a)=S$ for all $a\in S$.
However, we consider the members of $V$ to be marked by the letter
$a$, so that $V(a)$ is an unambiguous member of $V$ (strictly,
we are using the pairs $\{a,V(a)\},$ $(a\in S)$).

\textbf{Proposition 1.2 }For a finite semigroup $S$ the following
are equivalent:

(i) $S$ has a permutation matching;

(ii) $S$ is a transversal of $V=\{V(a)\}_{a\in S}$;

(iii) $|A|\leq|V(A)|$ for all $A\subseteq S$. 

\emph{Proof }(i) $\Rightarrow$(ii) Let $f:S\rightarrow S$ be a permutation
matching. Replacing $a$ by $f^{-1}(a)$ in the condition that $f(a)\in V(a)$
gives that $f^{-1}:S\rightarrow S$ is a permutation and $a\in V(f^{-1}(a))$
for all $a\in S$, which in turn is equivalent to the condition that
the function $g:S\rightarrow V$ whereby $g(a)=V(f^{-1}(a))$ is a
bijection such that $a\in g(a)$ for all $a\in S$ as 
\[
a\in g(a)\,\forall a\in S\Leftrightarrow a\in V(f^{-1}(a))\,\forall a\in S\Leftrightarrow f(a)\in V(a)\,\forall a\in S.
\]
 Hence if $f$ is a permutation matching then $g$ is a required bijection
that shows that $S$ is a transversal of $V$. 

(ii) $\Rightarrow$ (i) Suppose that $V$ has a transversal $\tau$
so there exists a bijection $g:\tau\rightarrow V$ such that $t\in g(t)$
for all $t\in\tau$. Since $\tau\subseteq S$ and $|S|=|V|$ it follows
that $\tau=S$ and so $g:S\rightarrow V$ and $a\in g(a)$ for all
$a\in S$. Define $f:S\rightarrow S$ by $f(a)=g^{-1}(V(a))$, which
is then a permutation on $S$. Moreover $f(a)\in g(f(a))=V(a)$ and
so $f(a)$ is a permutation matching of $S$.

(ii) $\Leftrightarrow$ (iii) As explained above, for any transversal
$\tau$ of $V$ we necessarily have $\tau=S.$ By Hall's Marriage
Lemma a transversal exists if and only if for every subset $A\subseteq\bigcup V(a)=S$,
the condition $|A|\leq|V(A)|$ is satisfied. $\qed$

\textbf{Example 1.3 }A finite regular semigroup with no permutation
matching. Consider the $7$-element aperiodic completely $0$-simple
semigroup $S=\{(i,j):1\leq i\leq2,\,1\leq j\leq3\}\cup\{0\}$, where
$E(S)=\{(1,2),(1,3),(2,1)\}\cup\{0\}$. Then $V\{(2,2),(2,3)\}=\{(1,1)\}$
and so Hall's Condition is violated on the set $A=\{(2,2),\,(2,3)\}$,
whence, by Proposition 1.2, $S$ has no permutation matching. This
example is minimal in view of the next result.

\textbf{Proposition 1.4 }Any regular semigroup with fewer than $7$
elements has an involution matching.

\emph{Proof }Suppose that $S$ is a finite regular semigroup with
$|S|\leq6$. If $S$ had just one ${\cal D}$-class, then $S$ would
be completely simple and so, being a union of groups, would have an
involution matching. Suppose then that $S$ is not a union of groups
and so has more than one ${\cal D}$-class. Then since $S$ is regular
but not a union of groups, it follows that $S$ possesses one ${\cal D}$-class
$D$ with more than one ${\cal R}$-class and more than one ${\cal L}$-class.
It now follows that $|D|\leq5$ and since $|D|$ is not prime it follows
that $D$ is a $2\times2$ `eggbox' with trivial ${\cal H}$-classes.
Moreover either each ${\cal R}$-class and each ${\cal L}$-class
consists of a pair of elements exactly one of which is idempotent
or $D$ has exactly one non-idempotent. In the first case each element
in $D$ has a unique inverse. In the 3-idempotent case it is also
clear that the four members of $D$ may be matched in mutually inverse
pairs (with two idempotents self inverse). The remaining ${\cal D}$-classes,
whether they number $1$ or $2,$ either consist of idempotents or
constitute a 2-element group. In any event it follows that an involution
matching for $S$ can be found. Therefore if $|S|\leq6$ then $S$
has an involution matching. $\qed$ 

We denote the class of finite regular semigroups with a permutation
matching (respectively involution matching) by ${\cal M}$ (respectively
${\cal N}$).

\textbf{Proposition 1.5 }Each of the classes ${\cal M}$ and ${\cal N}$
is closed under the taking of direct products but not under the taking
of regular subsemigroups or homomorphic images.

\emph{Proof.} Let $S,T\in{\cal M}$. Let $S=S'=\{s':s\in S\}$ and
$T=T'=\{t':t\in T\}$ denote permutation matchings of $S$ and of
$T$ respectively. For each $(s,t)\in S\times T$ put $(s,t)'=(s',t')\in V((s,t))$
whereupon $(')$ is a permutation of $S\times T$ that maps each element
to one of its inverses. Hence ${\cal M}$ is closed under the taking
of direct products. We note that if the matchings of $S$ and of $T$
are both involution matchings, then so is the matching $(')$ for
$S\times T$, so the same conclusion holds for the class ${\cal N}$.

We show in Theorem 2.12 that $T_{n}$ has a permutation matching.
By Example 1.3 we know that finite regular semigroups without permutation
matchings exist and since any finite semigroup may be embedded in
some $T_{n}$, it follows that ${\cal M}$ is not closed under the
taking of regular subsemigroups. However, consider the following example.

Let $T$ denote the aperiodic completely $0$-simple semigroup $T=\{(i,j):1\leq i,j\leq3\}\cup\{0\}$,
where $E(T)=\{(1,2),(1,3),(2,2),(2,3),(3,1)\}\cup\{0\}$. Note that
for each $a\in\{(1,1),(2,1)\}$ we have $V(a)=\{(3,2),(3,3)\}$. In
particular it follows that $T$ has an involution matching with pairs
$(1,1)\mapsto(3,3)$, $(2,1)\mapsto(3,2)$ and with $e\mapsto e$
for each $e\in E(T)$. 

We now note that the semigroup $S$ of Example 1.3, which has no permutation
matching, is a retract of $T$, meaning that $S$ is a subsemigroup
of $T$ and there is a homomorphism of $T$ onto $S$ that fixes each
member of $S$. The subsemigroup of $T$ that we identify with $S$
consists of the 7 members $S=\{(2,1),(2,2),(2,3),(3,1),(3,2),(3,3),0\}$.
Note that the relation $\rho=\{((1,j),(2,j)):\,1\leq j\;\leq3\}$
is a congruence on $T$ because a pair $(1,j)\in E(T)$ if and only
if $(2,j)\in E(T)$. A required retraction of $T$ onto $S$ is then
defined by $(1,j)\mapsto(2,j)\,(1\leq j\leq3)$, with each member
of $S$ fixed. Therefore $T\in{\cal N}$ but $S\not\in{\cal M}$ and
$S$ is both a regular subsemigroup of and a homomorphic image of
$T$. This completes the proof. $\qed$

\textbf{Theorem 1.6 }The following are equivalent for a finite semigroup
$S$:

(i) $S$ has a permutation matching;

(ii) $S$ has a permutation matching that preserves the ${\cal H}$-relation
(meaning that $\alpha{\cal H}\beta\Rightarrow\alpha'{\cal H}\beta'$);

(iii) each principal factor $D_{a}\cup\{0\}$ $(a\in S)$ has a permutation
matching;

(iv) each $0$-rectangular band $B=D_{a}\cup\{0\}/{\cal H}$ has a
permutation matching.

\emph{Proof} (i) $\Leftrightarrow$ (iii) Suppose that $S=S'=\{a',\cdots\}$
is a permutation matching of $S$ and let $D$ denote a ${\cal D}$-class
of $S$. For each $a\in D$ we have $V(a)\subseteq D$ so that $'|_{D}$
is a permutation of $D$; extending this by $0\mapsto0$ gives a permutation
matching for the principal factor $D\cup\{0\}$. Conversely, if each
principal factor has a permutation matching, the union of these matchings
over the set of ${\cal D}$-classes yields a permutation matching
of $S$.

(i) $\Rightarrow$ (iv) Let $\Lambda$ denote any collection of ${\cal H}$-classes
within $D=D_{a}$ and let $n$ be the common cardinal of the ${\cal H}$-classes
within $D$. For each $a\in H\in\Lambda$ we have that $H_{a'}\subseteq D\,(a'\in V(a))$.
Since $S$ has a permutation matching $'$, it follows that 
\[
|\{a\in\cup_{H\in\Lambda}H\}|=n|\Lambda|=|\{a':a\in\cup_{H\in\Lambda}H\}|\leq
\]
\begin{equation}
|\bigcup H_{a'}|\,(a\in\cup_{H\in\Lambda}H)=n|\{H_{a'}:a\in\cup_{H\in\Lambda}H\}|
\end{equation}

Let us denote the $0$-rectangular band $B$ by $(I\times J)\cup\{0\}$
and take any set $A\subseteq B$. The ${\cal H}$-classes within $D$
may be indexed by the members $(i,j)\in I\times J$. For each pair
of ${\cal H}$-classes $H_{i,j},H_{k,l}\subseteq D$ we have $H_{k,l}\subseteq V(H_{i,j})$
or $H_{k,l}\cap V(H_{i,j})=\emptyset$. For each $(i,j)\in A$ we
have that $(k,l)\in V((i,j))$ if and only if $H_{k,l}\subseteq V(H_{i,j})$.
Let $\Lambda=\{H_{i,j}:(i,j)\in A\}.$ Cancelling the common term
$n$ in (1) yields: 
\[
|A\setminus\{0\}|=|\Lambda|\leq|\{H_{a'}:a\in\cup_{h\in\Lambda}H\}|\leq|V(A\setminus\{0\})|;
\]
since $V(0)=\{0\}$ in $B$, it follows that $B$ satisfies Hall's
Condition and so $B$ has a permutation matching by Proposition 1.2. 

(iv) $\Rightarrow$ (ii) Let $D=D_{a}$ be the ${\cal D}$-class of
an arbitrary $a\in S$ consisting of ${\cal H}$-classes $H_{i,j}$
where the pairs $(i,j)$ form the set $I\times J$ say. Form a bipartite
graph $G$, the two defining independent sets of which are two copies,
$X$ and $Y$, of $I\times J$ with $(i,j)\in X$ adjacent to $(k,l)\in Y$
exactly when $(k,l)\in V((i,j)).$ By Proposition 1.2 applied to the
$0$-rectangular band $B$, each subset $A$ of $X$ is collectively
adjacent to at least $|A|$ vertices of $Y.$ By the Marriage Lemma,
it follows that there is a (one-to-one) matching of $X$ into $Y$,
which then induces a permutation of the set of ${\cal H}$-classes
of $D$ to itself, $H\mapsto H'$ say in such a way that each $a\in H$
has a unique inverse $a'\in H'$. The mapping whereby $a\mapsto a'$
is then a bijection of $H$ onto $H'$. 

Now list the ${\cal H}$-classes of $D$ as $H_{1},H_{2},\cdots.$
For each $a_{1}\in H_{1}$ let $\{a_{1}'\}=V(a)\cap H_{1}'$ and in
this way form a list $\{a_{1}',a_{2}',\cdots,a_{n}'\}$ say of all
the members of $H_{1}'$. We next select some $H_{2}\not\neq H_{1}$
and repeat the previous process, thus extending to a list of inverses
of the members of $H_{1}'\cup H_{2}'$ that exhausts this union. Continue
until all members of $X$ have been chosen. The final list $D'=\{a_{1}',a'_{2},\cdots\}$
will then be a list of inverses of all the members of $D$, which
is to say that $D'=D.$ Carrying out this procedure with each ${\cal D}$-class
of $S$ therefore results in a list $S'$ of inverses of the members
of $S$ that comprises the whole of $S$, that is $S'=S$ and $a\mapsto a'$
is a permutation matching of $S$ which, by construction, preserves
the ${\cal H}$-relation.

(ii) $\Rightarrow$ (i) is immediate. $\qed$

\textbf{Remarks 1.7} By above, if $S$ has a permutation matching
then $S$ has a permutation matching that preserves Green's ${\cal H}$
relation (and necessarily preserves the ${\cal D}$ relation) but
preservation of either the ${\cal L}$ or ${\cal R}$ relation cannot
be guaranteed even for inverse semigroups. To see this take any inverse
semigroup that is \emph{not} a union of groups (for example, the $5$-element
combinatorial Brandt semigroup $S={\cal M}[1,2,2,I_{2}]$, where $I_{2}$
is the $2\times2$ identity matrix). Let $a$ be a member of $S$
that does \emph{not }lie in a subgroup and let $e=aa^{-1}$. We have
$a{\cal R}e$ but $a^{-1}{\cal L}e=e^{-1}$. If $a^{-1}{\cal R}e^{-1}$
we would have $a^{-1}{\cal H}e$, which is a contradiction as $a^{-1}$
is not a group element. Hence ${\cal R}$ (and by symmetry also ${\cal L}$)
is not preserved by the (unique) permutation matching of $S$.

We may ask whether Proposition 1.6 goes through if we replace `permutation
matching' by `involution matching' throughout. In this case we do
see the implications ((i) $\Leftrightarrow$ (iii)) $\Leftarrow$
((ii) $\Leftrightarrow$ (iv)) but the missing forward implication
is not clear: to prove that (ii) $\Rightarrow$(iv) we note that we
may regard the non-zero members of $B$ as ${\cal H}$-classes of
some ${\cal D}$-class $D$ of $S$. We now have a well-defined involution
of $B$ in the mapping $H_{a}\mapsto H_{a'}$ and $0\mapsto0$, which
maps each member of $B$ to an inverse in $B$. Conversely, given
(iv), the involution on $B$ induces an involution by inverses between
pairs of ${\cal H}$-classes within $D$ so that (iv) implies (ii). 

We also do not have an example of a finite regular semigroup $S$
that has a permutation matching $f$ but no involution matching. In
any such example, all of the pemutation matchings $f$ must contain
an odd cycle that is free of idempotents, the reason being that the
set of members of any even cycle may be partitioned into mutually
inverse pairs while for any odd cycle that contains an idempotent
we may pair an idempotent in the cycle with itself, leaving an even
number of members that may be coupled into inverse pairs.

\textbf{Questions} Does $T_{n}$ have a permutation matching? An involution
matching? 

We shall give an affimative answer to the first question. A theorem
of Schein {[}21{]} (see {[}7, Theorem 6.2.4{]}) is that $T_{n}$ is
covered by its inverse subsemigroups. Indeed Schein proved that every
$a\in T_{n}$ has a \emph{strong inverse} $b\in V(a)$, which means
that the subsemigroup $\langle a,b\rangle$ is an inverse subsemigroup
of $T_{n}.$ However, there seems to be no guarantee that these subsemigroups
intersect in an inverse subsemigroup, which would allow us to infer
that $T_{n}$ does indeed have an involution matching.

\section{Graph Factor Interpretation}

All graphs $G=(V,E)$ under consideration will be finite graphs, perhaps
with loops (but not multiple edges). Recall that a $k$-factor of
$G$ is a subgraph $H=(V,E')$ of $G$ that is regular of degree $k$.
The most important cases are the $1$-factors and $2$-factors. A
$1$-factor $H$ is also known as a \emph{perfect matching }and is
essentially a set of disjoint edges of $G$ that span $G$, which
is to say covers all the vertices of $G$ so that $V(H)=V(G)$. A
$2$-factor on the other hand is a set of disjoint cycles of $G$
that collectively contain all of the graph's vertices. Since we allow
loops in our graphs, we shall here define a $1$-factor to be a subgraph
of $G$ that is a disjoint set of edges that cover all the vertices.
This allows for some of the edges to be loops (which correspond to
vertices of degree 2). To allow for $2$-cycles in our discussion,
we shall be considering $1,2$-factors of $G$, which are spanning
subgraphs $H$ of $G$ whose components are either edges or cycles.

\textbf{Definition 2.1} Let $S$ be a semigroup and let $G=G(S)=(S,E)$
be the \emph{graph of inverses}, which has vertex set $S$ and $uv$
is an (undirected) edge if $(u,v)\in V(S)$.

\textbf{Remarks 2.2 }A graph of inverses $G(S)$ may contain loops:
indeed there is a loop at $a\in G(S)$ if and only if $a=a^{3}$.
The semigroup $S$ is regular if and only if $G(S)$ has no isolated
vertices (vertices of degree $0$) and $S$ is an inverse semigroup
if and only if each component of $G(S)$ is a single edge (which may
be a loop: indeed in general there is a loop at every idempotent).
In general it is simple to show that the set of idempotents $E(S)$
of $S$ forms a semilattice if and only if Reg$(S)$ forms an inverse
subsemigroup of $S$. From this it follows that the components of
$G(S)$ each consist of at most $2$ vertices if and only if $E(S)$
forms a semilattice. 

This following type of graph was first introduced by Graham {[}5{]}
and later independently by Houghton {[}18{]} in relation to finite
$0$-simple semigroups although those authors explicitly worked in
terms of the Rees matrix representation of such semigroups. Graham
specifically addressed problems concerning maximal nilpotent subsemigroups
and maximal idempotent-generated subsemigroups while Houghton's paper
is based on the cohomology set of the graph $\Gamma$ over the underlying
group $G$ of the Rees matrix semigroup. A modern treatment of Graham's
approach appears in the book of Rhodes and Steinberg {[}20, Section
4.13.2{]} where the idempotent-generated subsemigroup of a Rees matrix
semigroup is studied in terms of the topology of this, its \emph{incidence
graph }in order to reveal results about varieties of finite semigroups
and general theorems such as that of Fitzgerald that the idempotent-generated
subsemigroup of a regular semigroup is itself regular {[}4{]}.

\textbf{Definition 2.3 }Let $S$ be a finite regular semigroup $S$. 

(a) The \emph{incidence graph }$GH(S)$ is a bipartite graph with
independent sets $\mathbf{L}$ and $\mathbf{R}$ of ${\cal L}$- and
${\cal R}$-classes respectively of $S$. An edge runs between $L\in\mathbf{L}$
and $R\in\mathbf{R}$ exactly when $L\cap R$ is a group.

(b) A ${\cal D}$-component of $GH(S)$ is the subgraph induced by
the set of all ${\cal L}$- and ${\cal R}$-classes in some ${\cal D}$-class
$D$ of $S$.

Since $S$ is regular, the set of ${\cal L}$-classes within a ${\cal D}$-class
$D$ of $S$ is collectively adjacent in $GH(S)$ to the set of ${\cal R}$-classes
within $D$ and vice-versa, so Definition 2.3(b) is symmetric in ${\cal L}$
and ${\cal R}$. It is not the case that a ${\cal D}$-component of
$GH(S)$ is necessarily connected but a ${\cal D}$-component of $GH(S)$
is a union of components of $GH(S)$. 

We shall call a ${\cal D}$-class of $S$ \emph{square }if it comprises
an equal number of ${\cal L}$- and ${\cal R}$-classes and we shall
call the semigroup $S$ itself \emph{square }if all of its ${\cal D}$-classes
are square.

\textbf{Theorem 2.4 }A ${\cal D}$-component subgraph $C$ of the
incidence graph $GH(S)$ of a finite regular semigroup $S$ has a
perfect matching if and only if the corresponding ${\cal D}$-class
$D$ is square and there is an ${\cal H}$-class preserving involution
matching of the principal factor $D\cup\{0\}$. Overall the graph
$GH(S)$ has a perfect matching if and only if $S$ is square and
$S$ has an ${\cal H}$-class preserving involution matching. 

\textbf{Remarks 2.5} Neither of the latter two conditions of the first
statement of the theorem imply the other. For example, the identity
mapping is an ${\cal H}$-class preserving involution matching for
any non-trivial right zero semigroup $S$, the unique ${\cal D}$-class
of which is then not square. Next, take the the $10$-element combinatorial
Rees matrix semigroup $S$ the non-zero ${\cal D}$-class of which
is the $3\times3$ array indexed by $\{1,2,3\}\times\{1,2,3\}$ whose
idempotents lie in positions $(1,j)$ and $(i,1)$ $(i,j\in\{1,2,3\})$.
Then $S$ is square but $V\{(2,2),(2,3),(3,2),(3,3)\}=\{(1,1)\}$,
so that Hall's Condition on the family of sets $V(a)$ is violated
and therefore, by Proposition 1.2, $S$ has no permutation matching.
Note also that in the semigroup $S$ of Example 1.3 there is a matching
of the members of $\mathbf{R}$ into those of $\mathbf{L}$ (but not
conversely) and again $S$ has no permutation matching.

\emph{Proof of Theorem 2.4} We prove the first statement of the theorem
from which the second follows by the argument of Theorem 1.6(i) $\Leftrightarrow$
(iii). Suppose there exists a perfect matching for a ${\cal D}$-component
$C$ of $GH(S)$ corresponding to some ${\cal D}$-class $D$ of $S$
and let $\mathbf{L}$, $\mathbf{R}$, and $\mathbf{H}$ denote the
respective collections of ${\cal L}$-, ${\cal R}$-, and ${\cal H}$-classes
of $D$. We shall write a typical ${\cal H}$-class $H$ as $H=L\cap R$,
where $L$ and $R$ are respectively the ${\cal L}$- and the ${\cal R}$-class
of $D$ that contain $H$. Let $f:\mathbf{L}\rightarrow\mathbf{R}$
denote a perfect matching of the ${\cal D}$-component $C$ of $GH(S)$
corresponding to $D$. It follows at once that $|\mathbf{L|=|R|}$
and so $D$ is square. 

Now introduce the mapping $':\mathbf{H}\rightarrow\mathbf{H}$ by
$H'=(L\cap R)'=Rf^{-1}\cap Lf$. By definition of the graph $GH(S)$,
we have that $L\cap Lf$ and $Rf^{-1}\cap R$ are group ${\cal H}$-classes
and so $H'$ is a set of inverses of $H$. (Indeed there is a perfect
matching between $H$ and $H'$ in the graph of inverses $G(S)$,
although $H=H'$ is possible in which case some of the edges of this
matching may be loops.) Note that the mapping $'$ is an involution
(and in particular a bijection) on $\mathbf{H}$ for we have:
\[
(H')'=((L\cap R)')'=(Rf^{-1}\cap Lf)'=(Lf)f^{-1}\cap(Rf^{-1})f=L\cap R=H.
\]
 Therefore we may define an ${\cal H}$-class preserving involution
matching, which we shall also denote by $'$, on the principal factor
of $D$ by taking an arbitrary ${\cal H}$-class $H\in\mathbf{H}$
and an arbitrary member $h\in H$ and defining $h\mapsto h'$ where
$h'\in H'$ is the unique inverse of $h$ that is to be found in $H'.$
This completes the proof of the forward implication of the first statement
of the theorem.

Conversely suppose that the principal factor of the ${\cal D}$-class
$D$ of $S$ is square so that $\mathbf{|L|=|R|}=n$ say, where $\mathbf{L}$
and $\mathbf{R}$ once again denote the respective collections of
${\cal L}$- and of ${\cal R}$-classes contained in $D$, which then
forms an $n\times n$ array of the ${\cal H}$-classes that comprise
the collection $\mathbf{H}$. Let us suppose further that the principal
factor of $D$ has an ${\cal H}$-class preserving involution matching
$':D\cup\{0\}\rightarrow D\cup\{0\}$. Let $T$ be any collection
of $k$ say members of $\mathbf{L}.$ The the union of the members
of $T$ is the union of $nk$ ${\cal H}$-classes of $D$, which is
mapped by the bijection $'$ onto a set of $nk$ ${\cal H}$-classes
of $D$. Since each member of $\mathbf{R}$ contains exactly $n$
${\cal H}$-classes, it follows that $T'$ meets at least $k$ of
the ${\cal R}$-classes of $D$. Now for each such ${\cal R}$-class
$R_{1}$ that meets $T'$ there exists some ${\cal H}$-class $H_{1}$
contained in some ${\cal L}$-class $L_{1}\subseteq T$ such that
$H_{1}'\subseteq R_{1}.$ From this it follows that $L_{1}\cap R_{1}$
is a group so that $L_{1}R_{1}$ is an edge in $GH(S)$ whence each
such $R_{1}$ is adjacent to some member of $T$. Hence the set of
$k$ ${\cal L}$-classes that form $T$ is collectively adjacent to
at least $k$ ${\cal R}$classes of $D$. It follows by Hall's Lemma
that there is a matching of $\mathbf{L}$ into $\mathbf{R}$ in the
${\cal D}$-component $C$ of $GH(S)$ corresponding to $D$, which,
since $|\mathbf{L|=|R|}$, is a perfect matching of the independent
sets $\mathbf{L}$ and $\mathbf{R}$ within each ${\cal D}$-class
$D$. $\qed$

\textbf{Examples 2.6 }An example that illustrates Theorem 2.4 is provided
by the semigroup of all orientation-preserving mappings $OP_{n}$
on an $n$-cycle $c=(0\,1\,\cdots\,n-1)$, which are the members $\alpha\in T_{n}$
such that the sequence $(0\alpha,1\alpha,\cdots,(n-1)\alpha)$ is
cyclic. As shown in {[}1{]}, $OP_{n}$ is a finite regular monoid
and the kernel classes of $\alpha\in OP_{n}$ are convex (that is
form intervals of the underlying cycle). Hence the set of kernel classes
(and hence the ${\cal R}$-class) of a mapping $\alpha\in OP_{n}$
of rank $k\geq2$ is determined by a strictly ascending list $P$
of $k$ integers $P=(p_{0}<p_{1}<\cdots<p_{k-1})$; the $p_{i}$ are
drawn from the integer interval $[n]=\{0,1,\cdots,n-1\}$ and each
$p_{i}$ equals the initial member of a kernel class of $\alpha$
when read in cylic order. Since each ${\cal L}$-class of $OP_{n}$
can also be identified with the common image of its members, which
corresponds to a $k$-set $A$ drawn from $[n]$, it follows that
there is a perfect matching for each ${\cal D}$-component of $GH(OP_{n})$
corresponding to a ${\cal D}$-class of rank $k\geq2$, defined by
matching the ${\cal L}$-class $L\in\mathbf{L}$ with the ${\cal R}$-class
$R\in\mathbf{R}$ that arises by taking the common range $A$ of members
of $L$ and letting $A$ act as the set $P$ of initial members of
kernel classes of the kernel partition that defines $R$. 

By our theorem, it follows that the ${\cal D}$-class $D$ of $OP_{n}$
of rank $k\geq2$ is square (for $D$ we have $\mathbf{|L|=|R|=}\binom{n}{k}$)
and the corresponding principal factor $D\cup\{0\}$ has an ${\cal H}$-class
preserving involution matching $'$. We can identify this involution
explicitly as follows. Let $\alpha\in OP_{n}$ be a mapping of rank
$k\geq2$ and let $\alpha$ lie in the ${\cal H}$-class $L\cap R$
say with $A=\{a_{0}<a_{1}<\cdots<a_{k-1}\}$ representing the common
image set of members of $L$ and let $P=\{p_{0}<p_{1}<\cdots<p_{k-1}\}$
represent the common set of initial members of kernel classes of $R$.
Furthermore, list the kernel classes of $R$ as $K_{0},K_{1},\cdots,K_{k-1}$,
where $p_{i}$ is the initial member of $K_{i}$ $(0\leq i\leq k-1)$.
The mapping $\alpha\in H$ is then specified by the choice of an integer
$r\in[k]$ where $K_{i}\alpha=a_{i\cdot c^{r}}$, where $c=c_{k}$
is the cycle $(0\,1\,\cdots\,k-1)$. The canonical inverse $\alpha'$
of $\alpha$ given by the involution matching of our theorem is then
the mapping $\beta$ with ${\cal L}$-class determined by the set
$P$ and ${\cal R}$-class determined by $A$ (so that the roles of
the sets $A$ and $P$ are interchanged when passing form $\alpha$
to $\beta$) and with $r$ replaced by $k-r$. Hence if we identify
$\alpha\in OP_{n}$ with the triple $\alpha=(A,P,r)$ then $\beta=(P,A,k-r)$.
It is now clear that the mapping $\alpha\mapsto\beta$ is an ${\cal H}$-class
preserving involution on $D$ (as an ${\cal H}$-class is determined
by the pair $(A,P)$) and we may check that $\beta\in V(\alpha)$
as follows. Let $t\in[n]$ with $t\in K_{i}$ say so that the initial
member of $K_{i}$ is $p_{i}$. Then we obtain:
\[
t\alpha\beta\alpha=p_{i}\alpha\beta\alpha=a_{i\cdot c^{r}}\beta\alpha=p_{(i\cdot c^{r})\cdot c^{k-r}}\alpha=p_{i\cdot c^{k}}\alpha=p_{i}\alpha=t\alpha,
\]
so that $\alpha\beta\alpha=\alpha$ and in the same way we have $\beta=\beta\alpha\beta$,
whence $\beta\in V(\alpha)$ and we may denote $\beta$ by $\alpha'$.
We conclude that $\alpha\mapsto\alpha'$ is an ${\cal H}$-class preserving
involution matching on $D$. 

The ${\cal D}$-class of all mappings of rank $1$ in $OP_{n}$ consists
of idempotents (indeed it comprises the set of all right zeros of
$OP_{n}$), which has an ${\cal H}$-class preserving involution matching
in the identity mapping. By taking the union of these matching across
all ${\cal D}$-classes of $OP_{n}$, we obtain an ${\cal H}$-class
preserving involution matching of $OP_{n}$. 

An example of a different kind is provided by the full linear monoid
$M_{n}(F)$ of all $n\times n$ matrices over a finite field $F$.
This monoid also satisfies the conditions of Theorem 2.4, and thus
also admits an involution matching.

\textbf{Proposition 2.7 }The semigroup $S$ has an involution matching
(respectively a permutation matching) if and only if $G=G(S)$ has
a $1$-factor (respectively a $1,2$-factor).

\emph{Proof}. Let $':S\rightarrow S$ be an involution matching of
$S$. Then the set $M=\{aa':a\in S\}\subseteq E(G)$. We then see
that the subgraph $H$ of $G$ defined by the set of edges $M$ is
a $1$-factor of $G$ as follows. Since $'$ has domain $S$, then
$V(H)=S.$ Any walk of length $2$ in $H$ has the form $a\rightarrow a'\rightarrow(a')'$
but since $(a')'=a$, this walk involves only a single edge, whence
$M$ comprises a set of disjoint edges. Conversely, let $H$ be a
$1$-factor of $G(S)$. Then each vertex $a\in G$ lies on a unique
edge $aa'$ of $G$ (with $a'=a$ possible) whence, by definition
of $G(S)$, $a\mapsto a'$ is an involution matching of of $S$.

Next suppose that $S$ has a permutation matching $':S\rightarrow S$,
which can then be regarded as a set of disjoint (oriented) cycles,
edges, and fixed points, the vertices of which cover $G(S)$. The
underlying set of loops (for $1$-cycles), edges (for $2$-cycles)
and cycles then represents a $1,2$-factor of $G(S)$. Conversely,
suppose that $H$ is a $1,2$-factor of $G(S)$. Choose any orientation
for each of the cycles of $H$ and treat each edge component of $H$
as a $2$-cycle, thereby defining a permutation $a\mapsto a'$ of
$S$ such that $a'\in V(a)$, which is by definition a permutation
matching of $S$. $\qed$

The main aim of this section is to show that the finite full transformation
semigoup $S=T_{n}$ on the base set $X_{n}=\{1,2,\cdots,n\}$ has
a permutation matching, which by Proposition 2.7 is the same as saying
that the graph of inverses, $G(S)$ has a $1,2$-factor. We reformulate
the notion of $1,2$-factor in terms of a certain bipartite graph
$G'$ known as the \emph{bipartite double cover }of $G$, also known
as the \emph{canoncial cover }or sometimes the \emph{Kronecker cover}
as $G'$ is realised as a certain product of $G$ with the single
edge $K_{2}$. The original source of this construction seems to be
the paper {[}2{]}.

\textbf{Proposition 2.8 }Let $G=(V,E)$ be a graph (with loops) and
let $V'$ be a disjoint copy of the vertex set $V$. Let $G'$ be
the bipartite graph with independent sets $V$ and $V'$ with $uv'\in E(G')$
if and only if $uv\in E(G)$. Then $G$ has a $1,2$-factor if and
only if $G'$ has a $1$-factor.

\emph{Proof} Suppose that $G'$ has a $1$-factor $H$. For each $u\in V(G)$
let $v'\in V(G')$ be such that $uv'\in E(H)$. Beginning with an
arbitrary vertex $u=u_{0}\in V(G)$ we may form a path $u_{0}\rightarrow u_{1}\rightarrow u_{2}\rightarrow\cdots$
in $G$ by defining for $i\geq1$, $u_{i}=v$ where $u_{i-1}v'$ is
the unique edge from $u_{i-1}$ in $H$. There is then a least value
of $i\geq0$ such that for some $j\geq1$ we have $u_{i}=u_{i+j}=v$,
say. If $i\geq1$ we would have that $u_{i-1}v'$ and $u_{i+j-1}v'$
were distinct edges of $H$ with a common vertex, contrary to $H$
being a $1$-factor of $G'$. Hence $u_{i+j}=u_{0}$, giving either
a cycle $C$ in $G$ containing $u_{0}$ or a walk of the form $u_{0}u_{1}u_{0}$.
Suppose recursively that we have constructed a set of cycles and edges
of $G$: $C_{1},C_{2},\cdots$ that are pairwise disjoint. Suppose
there remains a vertex $x=x_{0}$ of $V(G)$ that is not a vertex
of any of the $C_{i}$. Form a new cycle or edge in $G$: $x_{0}\rightarrow x_{1}\rightarrow\cdots\rightarrow x_{t}$$(t\geq0)$
as above. Suppose that some vertex $x_{i}$$(1\leq i)$ is such that
$x_{i}=u_{j}$, where $u_{j}\in V(C)$ where $C$ is one of the cycles
or edges already in our list of cycles of $G$. Take $i$ to be the
least index for which this is true. Then in $G'$ we have $x_{i-1}x_{i}'$
and $u_{j-1}u_{j}'$ are both edges of $H$ but $x_{i}'=u_{j}'$ while
$x_{i-1}\neq u_{j-1}$, contradicting that $H$ is a $1$-factor of
$G'$. Hence this eventuality does not arise and the process halts
with a required list of pairwise disjoint cycles and edges of $G$
that collectively span $G$.

Conversely, suppose that $G$ has a $1,2$-factor $H$, which then
comprises a set of pairwise disjoint cycles and edges $C_{1},C_{2},\cdots$
of $G$ that cover all of $V(G)$. Choose an orientation for each
cycle $C_{i}$, which then induces a mapping $':V(G)\rightarrow V(G)$
whereby $u\mapsto v'$, where $v$ is the unique vertex that follows
$u$ in the oriented cycle that contains the vertex $u$, or $v$
is the other vertex of the edge if $C_{i}$ is just an edge. The corresponding
set of edges $uv'$ of $G'$ then comprise a $1$-factor of $G'$.
$\qed$ 

\textbf{Lemma 2.9 }Let $G$ be a bipartite graph defined on two non-empty
disjoint independent sets $X$ and $Y$. Suppose that $X=\bigcup X_{i}$
and $Y=\bigcup Y_{i}$, where both unions are disjoint and comprise
an equal number of $t\geq1$ sets. Define the bipartite subgraphs
$G_{i}$ of $G$ by saying that $G_{i}$ has the pair of independent
sets $(X_{i},Y_{i})$ where $x\in X_{i}$ is adjacent to $y\in Y_{i}$
if $xy$ is an edge of $G$. Then $G$ has a perfect matching if each
of the graphs $G_{i}$ has an $m=m_{i}$-factor for some integer $m_{i}$. 

\emph{Proof }It suffices to check that each $G_{i}$ has a perfect
matching for, given this, the union of these perfect matchings will
be the required perfect matching for $G$. For a given $G_{i}$, let
$H$ be a regular subgraph of $G_{i}$ that spans $G_{i}$ so that
each vertex of $H$ has degree $m\geq1$ say. Take any subset of $k$
vertices of $X_{i},$ which are collectively adjacent to $l$ say
vertices in $Y_{i}$. Hence there are at least $km$ edges in $H$
incident with these $l$ vertices. Since every vertex in $H$ has
degree $m$, it follows that $lm\geq km\Rightarrow l\geq k$. It therefore
follows by Hall's Marriage Lemma that $G_{i}$ has a perfect matching.
(It follows also that each $|X_{i}|=|Y_{i}|$ and $|X|=|Y|$.) $\qed$

Recall the basic structure of Green's relations for $T_{n}$, which
are that $\alpha{\cal L\beta}$ if and only if $X\alpha=X\beta$,
$\alpha{\cal R}\beta$ if and only if ker $\alpha=$ ker $\beta$,
$\alpha{\cal D}\beta$ if and only if $|X\alpha|=|X\beta|$ with ${\cal D}={\cal J}$.
We shall write Ker $\alpha$ for the partition of $X$ induced by
ker $\alpha$. Since ${\cal H}={\cal R\cap{\cal L}}$ it follows that
an ${\cal H}$-class $H$ of $S=T_{n}$ is characterised by a kernel-range
pair $(\Pi,Y)$ where $Y\subseteq X$ and $\Pi$ is a partition of
$X$ with $|Y|=|\Pi|=k$ say, $(1\leq k\leq n)$; $k$ is the common
\emph{rank} of members in the ${\cal D}$-class containing $H$. 

It is convenient here however to introduce a new equivalence relation
${\cal {\cal Q}}$ on $T_{n}$ that lies between ${\cal R}$ and ${\cal D}$. 

\textbf{Definition 2.10 }For $\alpha\in T_{n}$ let us write Ker $\alpha=\{K_{1},K_{2},\cdots,K_{k}\}$
where the kernel classes of $\alpha$ are listed in ascending order
of cardinality. Let $P_{\alpha}=(p_{i})_{1\leq i\leq k}$, where $p_{i}=|K_{i}|$.
We shall say that $\alpha{\cal {\cal Q}}\beta$ in $T_{n}$ if $P_{\alpha}=P_{\beta}$.
We shall write $Q_{\alpha}$ for the ${\cal Q}$-class of $\alpha$.

Consider the double bipartite cover $G'$ of the graph of inverses
$G$ of $S=T_{n}$. For each ${\cal Q}$-class $Q_{\alpha}$ of $T_{n}$,
let $G_{\alpha}$ be the bipartite subgraph on the independent sets
$Q_{\alpha},Q_{\alpha}'$ where $uv'\in E(G_{i})$ if $u,v\in Q_{\alpha}$
and $uv\in G$. Our main result will follow once the next fact is
proved.

\textbf{Lemma 2.11} Each of the graphs $G_{\alpha}$ as given above
is $m$-regular for some positive integer $m=m_{\alpha}$ that depends
on $G_{\alpha}$. In particular, $G_{\alpha}$ is itself an $m$-factor
of itself.

\textbf{Theorem 2.12 }The semigroup $S=T_{n}$ has a permutation matching.

\emph{Proof} The ${\cal Q}$-classes of $S$ and their dashed counterparts
partition the independents sets $S$ and $S'$ that together comprise
the vertex set of $G'$ into an equal number of sets. The bipartite
graphs $G_{\alpha}$ based on the independent sets $Q_{\alpha}$ and
$Q_{\alpha}'$ defined above are then in accord with the hypotheses
of Lemma 2.9. Assuming the truth of Lemma 2.11, we may invoke Lemma
2.9 to conclude that $G'$ has a $1$-factor, whence by Proposition
2.8, $G(S)$ has a $1,2$-factor. But then by Proposition 2.7, we
conclude that $T_{n}$ has a permutation matching. $\qed$

\textbf{Lemma 2.13 }Let $H_{1},H_{2}$ be any two ${\cal H}$-classes
of $T_{n}$ defined by the respective kernel-range pairs $(\Pi_{1},Y_{1})$
and $(\Pi_{2},Y_{2})$. Then each $\alpha\in H_{1}$ has a (unique)
inverse $\beta\in H_{2}$ (and vice versa) if and only if $Y_{1}\,t\,\Pi_{2}$
and $Y_{2}\,t\,\Pi_{1}$. 

\emph{Proof }Quite generally by standard semigroup theory (consequences
of Green's Lemma), each member of $H_{1}$ will be an inverse of some
unique member of $H_{2}$ and vice versa if and only if the ${\cal H}$-classes
defined by the kernel-range pairs $(\Pi_{2},Y_{1})$ and $(\Pi_{1},Y_{2})$
are groups, which in turn occurs if and only if the given transversal
conditions are satisfied. $\qed$

\emph{Proof of Lemma 2.11}. Let $\alpha\in T_{n}$. We need to show
that the cardinal of $I=\{\beta\in T_{n}:\beta\in V(\alpha)\cap Q_{\alpha}\}$
is a positive integer $m$ that is independent of the choice of representative
of $Q_{\alpha}$. By Lemma 2.13, we may construct all members of $I$
as follows. Since we require that $\beta\in Q_{\alpha}$ we begin
with a list of `boxes' (named sets, the members of which are yet to
be specified) $L_{1},L_{2},\cdots,L_{k}$ where $k=|X\alpha|$, which
are to become the list of kernel classes of our mapping $\beta$.
All members of $Q_{\alpha}$ such that $X\alpha$ is a transversal
of the sets $L_{i}$ are formed by assigning exactly one member of
$X\alpha$ to each of the sets $L_{i}$ and then assigning all members
of $X\setminus X\alpha$ to the sets $L_{i}$ so that $|L_{i}|=p_{i}$.
Let us denote the number of ways in which this can be done by $l$,
noting that $l\geq1$. We are not asserting that the two stages of
assigning the members of $X\alpha$ to the $L_{i}$ and then assigning
the remaining members of $X$ to the $L_{i}$ are combinatorially
independent: indeed completion of the two stages will sometimes fail
to yield unique outcomes exactly when there is repetition among the
integers $p_{i}$. However, the value of $l$ depends only on $Q_{\alpha}$
and not on the choice of $\alpha$ as $l$ depends only on the cardinals
of the sets $L_{i}$ and the integer $k$ and not on the membership
of $X\alpha$. The number $l$ is then the number of partitions $\Pi$
of $X_{n}$ that correspond to kernel partitions of the members of
$Q_{\alpha}$ for which $X\alpha$ is a transversal of $\Pi$. 

All members of $I$ can now be formed as follows. Take any of the
$l$ partitions above to act as Ker $\beta$. To ensure that $\beta\in V(\alpha)$
we need to take $X\beta$ to be any transversal of Ker $\alpha$.
Again the number of choices $r$ available depends only on $P_{\alpha}$
and is independent of the the choice of representative of $Q_{\alpha}$;
indeed it is easy to see that $r=p_{1}p_{2}\cdots p_{k}\geq1$. The
required number $m=lr\geq1$ is then, by Lemma 2.13, the cardinal
of the set $I$. Therefore the subgraph $G_{i}$ in question is regular
of degree $m\geq1$, as required to complete the proof Lemma 2.11,
and hence also of Theorem 2.12. $\qed$

\textbf{Remarks 2.14} We note our proof shows additionally that $T_{n}$
has a permutation matching $\phi$ the restriction of which to each
${\cal Q}$-class of $T_{n}$ is also a permutation. We may also insist
that $\phi$ simultaneously preserves the ${\cal H}$-relation. To
see this, observe that for each ${\cal Q}$-class $Q$ in $T_{n}$,
the set $Q\cup\{0\}$ forms a subsemigroup (indeed a right ideal)
of the corresponding principal factor $D\cup\{0\}$. We may therefore
work through the proof of Theorem 1.6 with the semigroups $Q\cup\{0\}$
in place of the principal factors and ${\cal Q}$-preserving permutation
matchings throughout and so construct a permutation matching for $T_{n}$
that preserves both ${\cal Q}$- and ${\cal H}$-classes. 

We note that the semigroup generated by a fixed ${\cal Q}$-class
has been studied previously by Levi and her co-authors in several
papers: see for example {[}19{]}. 

A similar proof yields that the partial transformation semigroup $PT_{n}$,
which is isomorphic to the subsemigroup of all mappings that fix $0$
in $T_{X}$ where $X=\{0,1,2,\cdots,n\}$, also has a permutation
matching. We just need to run through the previous arguments, restricting
ourselves throughout to mappings that fix $0$. 

ACKNOWLEDGEMENT: The author would like to thank an anonymous referee
of this paper for some enlightening observations that led to its improvement.

\end{document}